\documentclass{article}
\usepackage{spconf,amsmath,graphicx}

\usepackage{graphicx}
\pdfoutput=1
\usepackage{amsfonts}
\usepackage{array}
\usepackage{amssymb}
\usepackage{amsmath}
\usepackage{algorithm}
\usepackage{caption}
\usepackage{subcaption}
\usepackage{algpseudocode}
\usepackage{bm}
\usepackage{color}
\usepackage[T1]{fontenc}
\usepackage{cite}
\newtheorem{thm}{Theorem}[section]

\newtheorem{defn}[thm]{Definition}
\newtheorem{cor}{Corollary}

\newcommand{\Exp}{\mathbb{E}}

\newcommand{\R}{\mathcal{R}}

\newcommand{\cG}{\mathcal{G}} 
\newcommand{\cV}{\mathcal{V}} 
\newcommand{\cE}{\mathcal{E}} 
\newcommand{\cL}{\mathcal{L}} 

\newcommand{\bA}{\mathbf{A}}

\newcommand{\bI}{\mathbf{I}}
\newcommand{\bL}{\mathbf{L}}

\newcommand{\bW}{\mathbf{W}}

\newcommand{\eqdef}{:=}

\newcommand{\cO}{{\cal O}}

\newcommand{\mA}{{\bf A}}

\newcommand{\mW}{{\bf W}}

\newcommand{\range}[1]{\operatorname{Range}\left(#1\right)}

\usepackage[colorinlistoftodos,bordercolor=orange,backgroundcolor=orange!20,linecolor=orange,textsize=scriptsize]{todonotes}

\title{Provably Accelerated Randomized Gossip Algorithms}
\name{Nicolas Loizou$^{a,b}$ \qquad Michael Rabbat${}^b$ \qquad Peter Richt\'{a}rik${}^{a, c}$}
\address{${}^c$ \em The University of Edinburgh, UK\\${}^b$ \em Facebook AI Research, Montreal \\${}^c$ \em KAUST, KSA\\}
\begin{document}
%
\maketitle
\begin{abstract}
In this work we present novel provably accelerated gossip algorithms for solving the average consensus problem. The proposed protocols are inspired from the recently developed accelerated variants of the randomized Kaczmarz method - a popular method for solving linear systems. In each gossip iteration all nodes of the network update their values but only a pair of them exchange their private information. Numerical experiments on popular wireless sensor networks showing the benefits of our protocols are also presented.
\end{abstract}
\section{Introduction}
\label{sec:intro}
Distributed averaging is a fundamental problem in the area of distributed computing and multi-agent systems~\cite{degroot1974reaching,tsitsiklis1986distributed}. Randomized gossip algorithms are one of the most popular class of methods for solving it. The seminal 2006 paper of Boyd et al.\ \cite{boyd2006randomized} on randomized gossip algorithms motivated a flurry of subsequent research, and now gossip algorithms appear in many applications, including distributed data fusion in sensor networks \cite{xiao2005scheme}, load balancing \cite{cybenko1989dynamic} and clock synchronization \cite{freris2012fast}.  The development and design of efficient gossip algorithms was studied extensively in the last decade. For a survey of selected relevant work prior to 2010, we refer the reader to the survey~\cite{dimakis2010gossip}. For more recent results on randomized gossip algorithms we suggest \cite{zouzias2015randomized, liu2013analysis,olshevsky2014linear,LoizouRichtarik, nedic2018network, aybat2017decentralized}. 
See also  \cite{dimakis2008geographic, aysal2009broadcast, olshevsky2009convergence,hanzely2017privacy}. 

In the literature of gossip algorithms, an important task is the design of fast and efficient algorithms. Surprisingly, to the best of our knowledge, there are no variants of gossip algorithms that converge to consensus with an accelerated linear rate. In this work, our focus is precisely this. We design two provably accelerated randomized gossip protocols which converge to consensus fast.

\vspace{0.2cm}
\noindent \textbf{The average consensus problem.}
In the average consensus (AC) problem we are given an undirected connected network $\cG=(\cV,\cE)$ with node set $\cV=\{1,2,\dots,n\}$ and edges $\cE$. Each node $i \in \cV$ holds a local value $c_i \in \R$. The goal of AC is for every node to compute the average of these private values, $\bar{c}\eqdef\tfrac{1}{n}\sum_i c_i$, in a distributed fashion. That is, the exchange of information can only occur between connected nodes (neighbors). 

\vspace{0.2cm}
\noindent \textbf{Main contributions.}
In this work, building upon a recent framework for the design and analysis of randomized gossip algorithms \cite{LoizouRichtarik, loizou2018accelerated}, we present two novel and provably accelerated randomized gossip protocols where in each step all nodes of the network update their values using their own information but only a pair of them exchange messages. The accelerated convergence rates of the proposed protocols are obtained by establishing a connection with the area of accelerated randomized Kaczmarz methods for solving consistent linear systems.

To the best of our knowledge, our protocols are the first randomized gossip algorithms that converge to consensus with an accelerated linear rate. The theoretical results are validated via computational testing on typical network topologies. 

\vspace{0.2cm}
\noindent \textbf{Structure of the paper.}
Section~\ref{sec:TechnicalPreliminaries} introduces important technical preliminaries and the necessary background for understanding of our methods. Two accelerated variants of the randomized Kaczmarz (RK) method for solving linear systems and their theoretical convergence results are described.  In Section~\ref{sec:AccGossip} we present the two provably accelerated gossip protocols, along with some remarks on their implementation. Numerical evaluation of the new gossip protocols is presented in Section~\ref{sec:NumericEval}. Finally, concluding remarks are given in Section~\ref{sec:Conclusion}.

\vspace{0.2cm}
\noindent \textbf{Notation.}
The following notational conventions are used in this paper. We write $[n]\eqdef \{1,2, \dots ,n\}$.  Boldface upper-case letters denote matrices; $\bI$ is the identity matrix. 
By $\cL$ we denote the solution set of the linear system $\bA x=b$, where $\bA \in \R^{m\times n}$ and $b\in \R^m$. By $\bA_{i:}$ and $\bA_{:j}$ we indicate the $i_{th}$ row and the $j_{th}$ column of matrix $\bA$, respectively. Throughout the paper, $x^*$ is the projection of $x^0$ onto $\cL$ (that is, $x^*$  is the solution of the best approximation problem; see equation \eqref{best approximation}). 
With $\lambda_{\min}^+(\cdot)$ we indicate the smallest nonzero eigenvalue of matrix $(\cdot)$. $\|\cdot \|$ and $\|\cdot \|_F$ are used to denote the Euclidean norm and  the Frobenius norm, respectively. Finally, $x^k = (x^k_1,\dots,x^k_n) \in \R^n$ represents the vector with the local values of the $n$ nodes of the network at the $k^{th}$ iteration. Here, $x_i^{k}$ denotes the value of node $i \in [n]$ at the $k^{th}$ iteration.

\section{Technical Preliminaries}
\label{sec:TechnicalPreliminaries}
In this section we present the connections between the randomized Kaczmarz methods for solving linear systems and the gossip algorithms for solving the AC problem, as discussed in more details in \cite{LoizouRichtarik, loizou2018accelerated}. In particular, we focus on the presentation of the two recently proposed accelerated variants of Kaczmarz methods and on their theoretical convergence analysis. 

\subsection{Kaczmarz-type methods and gossip algorithms}
Kaczmarz-type methods are popular algorithms for solving linear systems $\bA x =b$ with many equations. The randomized Kaczmarz method (RK) for solving consistent linear systems was first proposed and proved to converge with linear rate in \cite{RK}. This work triggered much research into developing and analyzing randomized linear solvers and several improved variants of RK have been proposed \cite{needell2010randomized, RBK, eldar2011acceleration, MaConvergence15, zouzias2013randomized, l2015randomized, schopfer2016linear, liu2016accelerated, loizou2017linearly}.

In particular, in its simplest version, RK works as follows; In each step, one row $\bA_{i:}$ of matrix $\bA$ is sampled with probability $p_i>0$ and then is used to obtain the next iterate by following the update rule:
\begin{equation}
\label{RK}
x^{k+1}=x^k - \tfrac{\bA_{i :} x^k -b_{i}}{\|\bA_{i :}\|_2^2} \bA_{i :}^ \top  .
\end{equation}
For the case of consistent linear systems, it was shown that RK and its variants solves the following problem (known as best approximation problem) \cite{gower2015randomized, gower2015stochastic, loizou2017momentum} :
\begin{equation}
\label{best approximation}
\min_{x = (x_1,\dots, x_n) \in \R^n} \tfrac{1}{2} \|x-x^0\|^2 
\quad \text{subject to}  \quad \bA x = b.
\end{equation}
where $x_0$ is the initial vector of the method.

In \cite{LoizouRichtarik} it was shown how RK works as a gossip algorithm when applied to a special linear system encoding the underlying network. The following definition is used to describe the class of linear systems considered here.
\begin{defn}[\cite{LoizouRichtarik}] 
\label{definition}
A linear system  $\bA x = b$ is called an ``average consensus (AC) system'' when  all solutions $x$ satisfy that $x_i = x_j$ for all $(i,j) \in \cE$.
\end{defn}
Many linear systems satisfy the above definition. In this work we focus on the case where $b=0$ and $\bA \in \R^{|\cE| \times n}$ is the incidence matrix of  $\cG$ (or its normalized form where $\|\bA_{i:}\|=1$). In this case, the row of the system  corresponding to edge $(i,j)$ directly encodes the constraint $x_i=x_j$. 

Since the right hand side of the above system is $b=0$, the update rule of equation \eqref{RK} simplifies to:
$x^{k+1}=x^k - \tfrac{\bA_{i :} x^k}{\|\bA_{i :}\|_2^2} \bA_{i :}^ \top = \left[\bI - \tfrac{\bA_{i :}^ \top  \bA_{i :}}{\|\bA_{i :}\|_2^2} \right] x^k.$ In the case that the starting point is $x_0=c$ it can be shown that RK solves the average consensus probem and that the above udpate rule is equivalent with the pairwise randomized gossip algorithm of \cite{boyd2006randomized} (see \cite{LoizouRichtarik} for more details). The convergence performance of RK for solving the best approximation problem (and as a result the average consensus problem) is described by the following theorem.
\begin{thm}[\cite{gower2015randomized,gower2015stochastic}]
\label{theoremRK}
Let $\{x^k\}$ be the iterates produced by \eqref{RK}. Then
$\Exp[\|x^k-x^*\|^2]\leq \rho^k \|x^0-x^*\|^2,$
where $\rho \eqdef 1 - \lambda_{\min}^+(\frac{\bA^\top \bA}{\|\bA\|^2_F}) \in [0,1]$. 
\end{thm}
\subsection{Accelerated Kaczmarz methods}
\label{AcceleratedVariants}
There are two different but very similar ways to accelerate the randomized Kaczmarz method. The first paper that proves asymptotic convergence with an accelerated linear rate is \cite{liu2016accelerated}. The proof technique is similar to the framework developed by Nesterov in \cite{nesterov2012efficiency} for the acceleration of coordinate descent methods.  In \cite{tu2017breaking,gower2018accelerated} a modified version for the selection of the parameters was proposed and a non-asymptotic accelerated linear rate was established. In Algorithm~\ref{alg:AccKaczmarz}, pseudocode of the Accelerated Kaczmarz method (AccRK) is presented where both variants can be cast as special cases, by choosing the parameters with the correct way. 
\begin{algorithm}[!h]
\begin{algorithmic}[1]
\State Data: Matrix $\mA\in \R^{m\times n}$; vector $b\in \R^m$
\State Choose $x^0\in \R^n$ and set $v^0 = x^0$
\State Parameters: 
Evaluate the sequences of the scalars $\alpha_k, \beta_k ,\gamma_k$ following one of two possible options.
\For {$k = 0, 1, 2, \dots, K$}
 \State $y^k = \alpha_k v^k + (1-\alpha_k) x^k$ 
\State Draw a fresh sample $i_k \in [m]$ with equal probability 
\State $x^{k+1} = y^k - \tfrac{\bA_{i_k :} y^k -b_{i_k}}{\|\bA_{i_k :}\|_2^2} \bA_{i_k :}^ \top$ 
\State $v^{k+1} = \beta_k v^k + (1-\beta_k) y^k -\gamma_k \tfrac{\bA_{i_k :} y^k -b_{i_k}}{\|\bA_{i_k :}\|_2^2} \bA_{i_k :}^ \top$
\EndFor
\end{algorithmic}
\caption{Accelerated Randomized Kaczmarz Method (AccRK)}
\label{alg:AccKaczmarz}
\end{algorithm}
There are two options for selecting the parameters, which we describe next.
\begin{enumerate}
\item From \cite{liu2016accelerated}: 
Choose $\lambda \in [0,\lambda_{\min}^+(\bA^\top \bA)]$ and set $\gamma_{-1}=0$.
Generate the sequence $\{\gamma_k: k=0,1,\dots, K+1\}$ by choosing $\gamma_k$ to be the largest root of $$\gamma_k^2-\frac{\gamma_k}{m}=(1-\frac{\gamma_k}\lambda{m})\gamma_{k-1}^2$$ and generate the sequences $\{\alpha	_k :  k=0,1,\dots,K+1\}$ and $\{\beta_k :  k=0,1,\dots,K+1\}$ by setting $$\alpha_k=\frac{m-\gamma_k\lambda}{\gamma_k(m^2-\lambda)}, \quad \beta_k=1-\frac{\gamma_k \lambda}{m}.$$
\item From \cite{gower2018accelerated}: Let 
\begin{equation}
\label{thenu}
\nu= \max_{u \in \range{\mA^\top}   } \frac{ u^\top \left[\sum_{i=1}^m  \mA_{i:}^\top \mA_{i:} (\mA^\top \mA)^\dagger \mA_{i:}^\top \mA_{i:} \right]u }{ u^\top \frac{\bA^\top \bA }{m}u }.
\end{equation}
Choose the three sequences to be fixed constants as follows:
$\beta_k =\beta = 1-\sqrt{ \frac{\lambda_{\min}^+(\bW)}{\nu}} $, \;$\gamma_k=\gamma = \sqrt{ \frac{1}{\lambda_{\min}^+(\bW) \nu}} $, \; $\alpha_k=\alpha =  \frac{1}{1+\gamma \nu} \in (0,1)$ where $\bW=\frac{\bA^\top \bA}{m}$.
\end{enumerate}
\subsection{Theoretical guarantees of AccRK}
The two variants (Option 1 and Option 2) of AccRK are closely related, however their convergence analyses are different. Below we present the theoretical guarantees of the two options as presented in \cite{liu2016accelerated} and \cite{gower2018accelerated}.
\begin{thm}[\cite{liu2016accelerated}]
Let $\{x^k\}_{k=0}^\infty$ be the sequence of random iterates produced by Algorithm~\ref{alg:AccKaczmarz} with the Option 1 for the parameters. Let $\lambda \in [0,\lambda_{\min}^+(\bA^\top \bA)]$ and define $\sigma_1=1+\frac{\sqrt{\lambda}}{2m}$ and $\sigma_2=1-\frac{\sqrt{\lambda}}{2m}$. Then for any $k\geq 0$ we have that:
$$\Exp[\|x^{k+1}-x^*\|^2 ] \leq \frac{4 \lambda}{(\sigma_1^{k+1}-\sigma^{k+1}_2)^2}\|x^0-x^*\|^2_{(\bA^\top \bA)^+}.$$
\end{thm}
\begin{cor}[\cite{liu2016accelerated}]
Note that as $k\rightarrow\infty$, we have that $\sigma^k_2\rightarrow 0$. This means that the decrease of the right hand side is governed mainly by the behavior of the term $\sigma_1$ in the denominator and as a result the method converge \emph{asymptotically} with a decrease factor per iteration:
$\sigma_1^{-2}=(1+\frac{\sqrt{\lambda}}{2m})^{-2}\approx 1-\frac{\sqrt{\lambda}}{m}.$
\end{cor}
Thus, by choosing $\lambda=\lambda_{\min}^+$ and for the case that $\lambda_{\min}^+$ is small, Algorithm~\ref{alg:AccKaczmarz} will have significantly faster convergence rate than RK. Note that the above convergence results hold for normalized matrices $\bA \in \R^{m \times n}$, that is matrices that have $\|\bA_{i:}\|=1$ for any $i \in m$. 
\begin{thm}[\cite{ gower2018accelerated}]
\label{theorem2}
Let $\bW=\frac{\bA^\top \bA}{m}$ and assume that ${\rm Null}(\mW) = {\rm Null}(\mA)$. Let $\{x^k,y^k,v^k\}$ be the iterates  of Algorithm~\ref{alg:AccKaczmarz} with the Option 2 for the parameters.  Then 
$$\Psi^k \leq \left(1-\sqrt{\lambda_{\min}^+(\bW)/\nu}\right)^k \Psi^0$$
where $\Psi^k =\Exp \left[\| v^k - x^*\|_{\mW^\dagger}^2+  \frac{1}{\mu  } \|x^{k}-x^*\|^2 \right]$
\end{thm}
The above result implies that Algorithm~\ref{alg:AccKaczmarz} converges linearly with rate $1-\sqrt{\lambda_{\min}^+(\bW)/\nu}$, which translates to a total of $O\left(\sqrt{\nu/\lambda_{\min}^+(\bW)}\log(1/\epsilon)\right)$ iterations to bring the quantity $\Psi^k$ below $\epsilon>0$. It can be shown that $1 \leq \nu \leq 1/\lambda_{\min}^+(\bW)$, (Lemma 2 in \cite{gower2018accelerated}) where $\nu$ is as defined in \eqref{thenu}. Thus,
$\sqrt{\frac{1}{\lambda_{\min}^+(\bW)}} \leq \sqrt{\frac{\nu}{\lambda_{\min}^+(\bW)}} \leq \frac{1}{\lambda_{\min}^+(\bW)}$
which means that the rate of AccRK (Option 2) is always better than that of the RK which (see Theorem~\ref{theoremRK}) is equal to $\cO(1/\lambda_{\min}^+(\bW)\log(1/\epsilon))$ for normalized matrices ($\|\bA\|^2_F=m$).
\section{Accelerated randomized gossip algorithms}
\label{sec:AccGossip}
In the previous section we presented the complexity analysis guarantees of AccRK for solving consistent linear systems with normalized matrices. Now, let us explain how the two options of AccRK behave as gossip algorithms when they are used to solve the linear system $\bA x=0$ where $\bA \in \R^{|\cE| \times n}$ is the normalized incidence matrix of the network. That is, each row $e=(i,j)$ of $\bA$ can be represented as $ (\bA_{e:})^\top= \frac{1}{\sqrt{2}}(e_i -e_j)$ where $e_i$ (resp.$e_j$) is the $i^{th}$ (resp. $j^{th}$) unit coordinate vector in $\R^{n}$.

By using this particular linear system, the expression $\tfrac{\bA_{i :} y^k -b_{i}}{\|\bA_{i :}\|_2^2} \bA_{i :}^ \top$ that appears in steps 8 and 9 of AccRK takes the following form when the row $e=(i,j) \in \cE$ is sampled:
$\tfrac{\bA_{e :} y^k -b_{i}}{\|\bA_{e :}\|_2^2} \bA_{e :}^ \top \overset{b=0}{=} \tfrac{\bA_{e :} y^k}{\|\bA_{e :}\|_2^2} \bA_{e :}^ \top \overset{\text{form of A}}{=} \frac{ y_i^k - y_j^k}{2}(e_i-e_j).$

Let $\bL$ be the Laplacian matrix of the network. For solving the above AC system (see Definition~\ref{definition}), the simple RK requires $O\left((\frac{2m}{\lambda_{\min}^+(\bL)})\log(1/\epsilon)\right)$ iterations to achieve expected accuracy $\epsilon>0$. To understand the acceleration in the gossip framework this should be compared to the $O(m\sqrt{2/\lambda_{\min}^+(\bL)} \log(1/\epsilon))$ of AccRK (Option 1) and the $O(\sqrt{2m \nu/\lambda_{\min}^+(\bL)} \log(1/\epsilon))$ of AccRK (Option 2).

Algorithm~\ref{alg:acceleratedNew} describes in a single framework how the two variants of AccRK of Section~\ref{AcceleratedVariants} behave as gossip algorithms when are used to solve the above linear system. Note that each node $\ell \in \cV$ of the network have two local registers to save the quantities $v^k_\ell$ and $x^k_\ell$. In each step using these two values every node $\ell \in \cV$ of the network (activated or not) computes the quantity $y^k_\ell =\alpha_k v^k_\ell + (1-\alpha_k) x^k_\ell$. Then in the $k^{th}$ iteration the activated nodes $i$ and $j$ of the randomly selected edge $e=(i,j)$ exchange their values $y^k_i$ and $y^k_j$ and update the values of $x^k_i$,  $x^k_j$ and $v^k_i$, $v^k_j$ as shown in Algorithm~\ref{alg:acceleratedNew}. The rest of the nodes use only their own $y^k_\ell$ to update the values of $v^k_i$ and $x^k_i$ without communicate with any other node.

The parameter $\lambda^+_{\min}(\bL)$ can be estimated by all nodes in a decentralized manner using the method described in~\cite{charalambous2016distributed}. In order to implement this algorithm, we assume that all nodes have synchronized clocks and that they know the rate at which gossip updates are performed, so that inactive nodes also update their local values. This may not be feasible in all applications, but when it is possible (e.g., if nodes are equipped with inexpensive GPS receivers, or have reliable clocks) then they can benefit from the significant speedup achieved.
\begin{algorithm}[!h]
\begin{algorithmic}[1]
\State Data: Matrix $\mA\in \R^{m\times n}$ be the normalized incidence matrix; vector $b=0\in \R^m$
\State Choose $x_0\in \R^n$ and set $v_0 = x_0$
\State Parameters: 
Evaluate the sequences of the scalars $\alpha_k, \beta_k ,\gamma_k$ following one of two possible options.
\For {$k = 0, 1, 2, \dots, K$}
\State Each node $\ell \in \cV$ evaluate $y^k_\ell = \alpha_k v^k_\ell + (1-\alpha_k) x^k_\ell$.
\State Pick an edge $e=(i,j)$ uniformly at random.
\State Then the nodes update their values as follows:
\begin{itemize}
\item The selected node $i$ and node $j$: 
$$x^{k+1}_i = x^{k+1}_j  = (y^k_i+y^k_j)/2$$ 
$$v^{k+1}_i = \beta_k v^k_i + (1-\beta_k) y^k_i -\gamma_k (y^k_i-y^k_j)/2$$
$$v^{k+1}_j = \beta_k v^k_j + (1-\beta_k) y^k_j -\gamma_k (y^k_j-y^k_i)/2$$
\item Any other node $\ell \in \cV$:
$$x^{k+1}_\ell = y^k_\ell \quad,\quad v^{k+1}_\ell = \beta_k v^k_\ell + (1-\beta_k) y^k_\ell$$
\end{itemize}
\EndFor
\end{algorithmic}
\caption{Accelerated Randomized Gossip Algorithm (AccGossip)}
\label{alg:acceleratedNew}
\end{algorithm}

\noindent \textbf{Related work on accelerated gossip algorithms:} The idea of having gossip updates in a network with two registers in each node is not new. It was first proposed in \cite{cao2006accelerated} and its analysis under strong conditions was presented in \cite{liu2013analysis}. There local memory is exploited by installing shift registers at each agent where the first register stores the agent's current value and the second the agent's value before the latest update. In \cite{loizou2018accelerated}, the Stochastic Heavy Ball (SHB) method is used for solving the AC problem and an accelerated method is proposed which was shown to be in practice faster than the algorithm of \cite{cao2006accelerated,liu2013analysis}. \cite{loizou2018accelerated} is the first paper that presents gossip algorithms where in each step all nodes of the network update their values but only a subset of them exchange their private values.
\section{Numerical Evaluation}
\label{sec:NumericEval}
We devote this section to numerically evaluate the performance of the proposed accelerated gossip protocols. 
In all of our experiments we compare the simple RK (equivalent to pairwise gossip algorithm of \cite{boyd2006randomized}) the Stochastic Heavy Ball method (SHB) proposed in \cite{loizou2018accelerated} and the AccRK (Algorithm~\ref{alg:acceleratedNew}) with the two options for the selection of the parameters presented in Section~\ref{AcceleratedVariants}. In comparing the methods we use the relative error measure $\|x^k-x^*\|^2 / \|x^0-x^*\|^2 $ where the starting vector of values $x^0=c$ is taken to be always Gaussian vector.  For all of our experiments the horizontal axis represents the number of iterations. The networks used in the experiments are the cycle (ring graph), the 2-dimension grid and the randomized geometric graph (RGG) with radius $r=\sqrt{\log(n)/n}$. Code was written in Julia 0.6.3. 

For the implementation of SHB we use the same parameters with the ones used in \cite{loizou2018accelerated}. For the AccRK (Option 1) we use $\lambda=\lambda_{\min}^+(\bA^\top\bA)$. Note that for all networks under study the two proposed protocols are faster than both the pairwise gossip algorithm of \cite{boyd2006randomized} and the SHB of \cite{loizou2018accelerated}.

\begin{figure}[t]
\centering
\begin{subfigure}{.24\textwidth}
  \centering
  \includegraphics[width=1\linewidth]{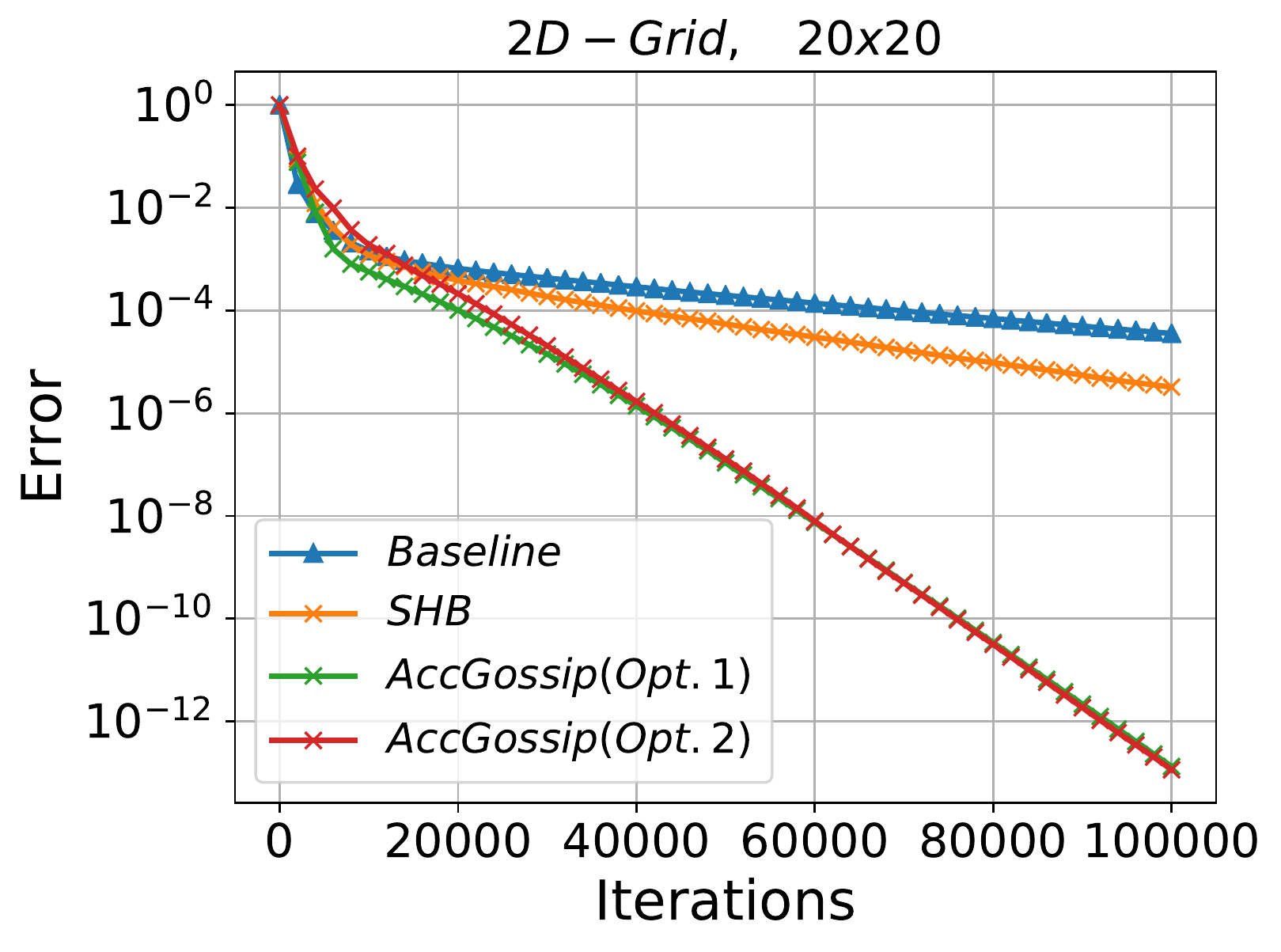}
\end{subfigure}%
\begin{subfigure}{.24\textwidth}
  \centering
  \includegraphics[width=1\linewidth]{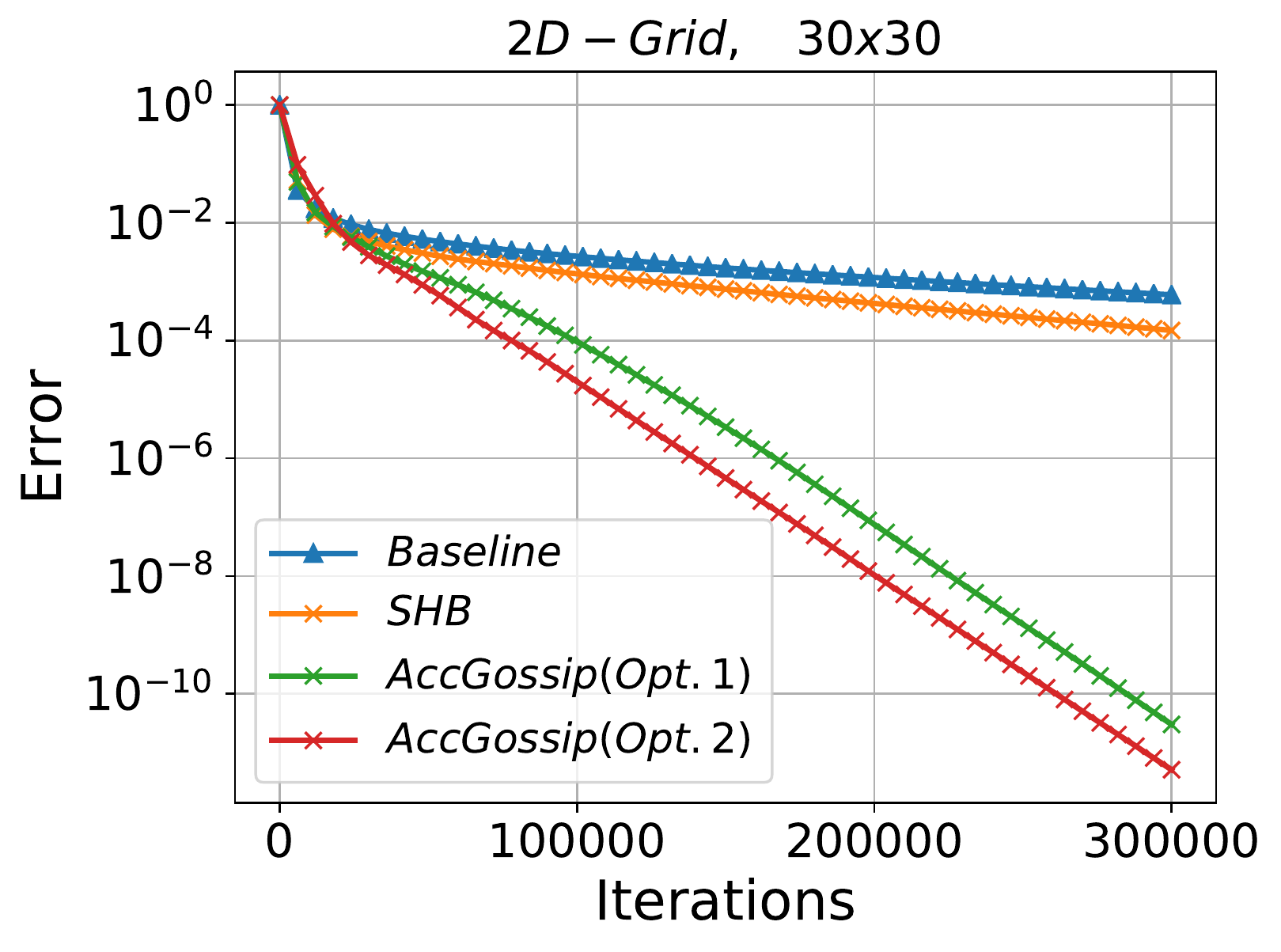}
\end{subfigure}\\
\begin{subfigure}{.24\textwidth}
  \centering
  \includegraphics[width=1\linewidth]{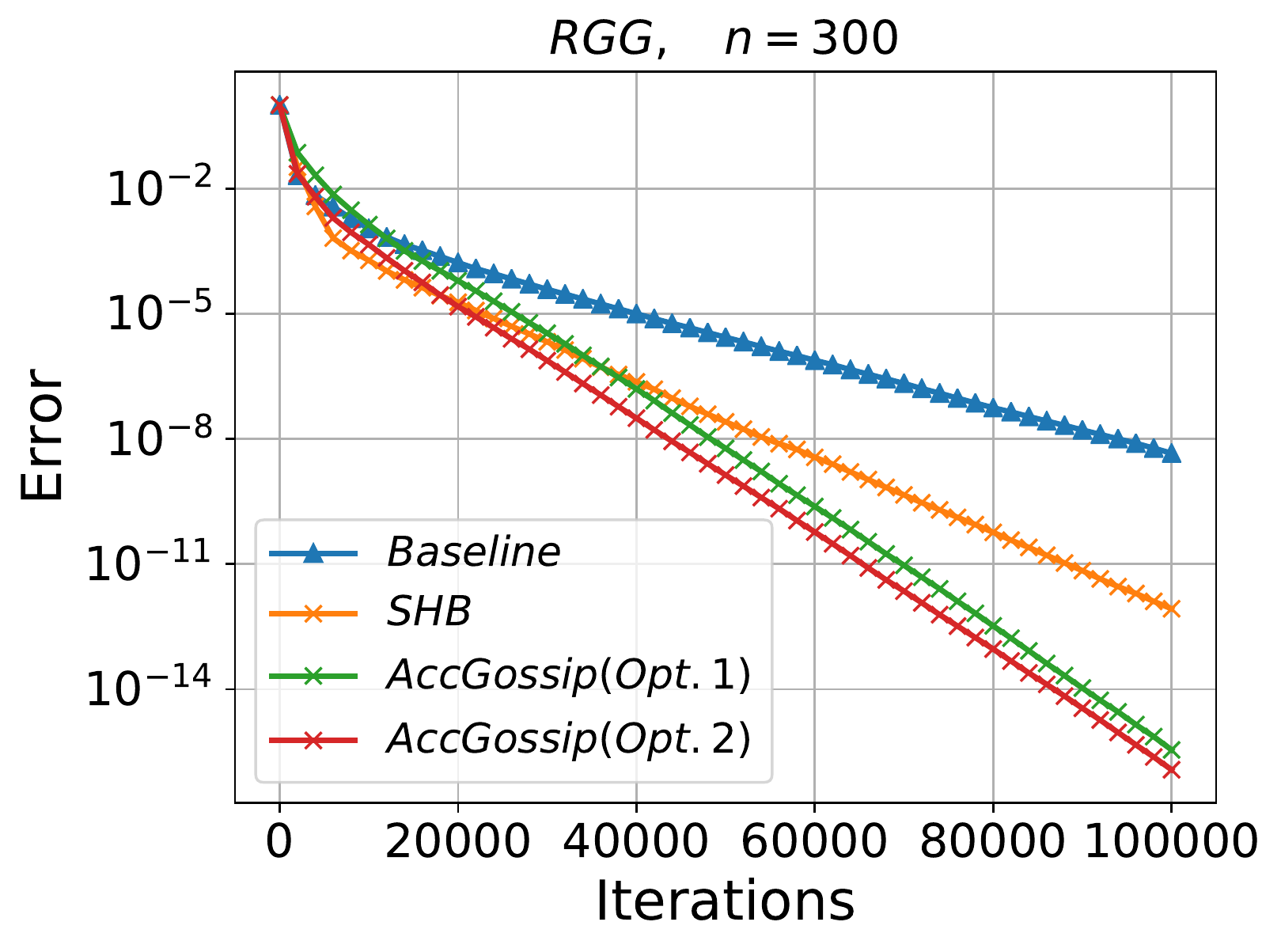}
\end{subfigure}%
\begin{subfigure}{.24\textwidth}
  \centering
  \includegraphics[width=1\linewidth]{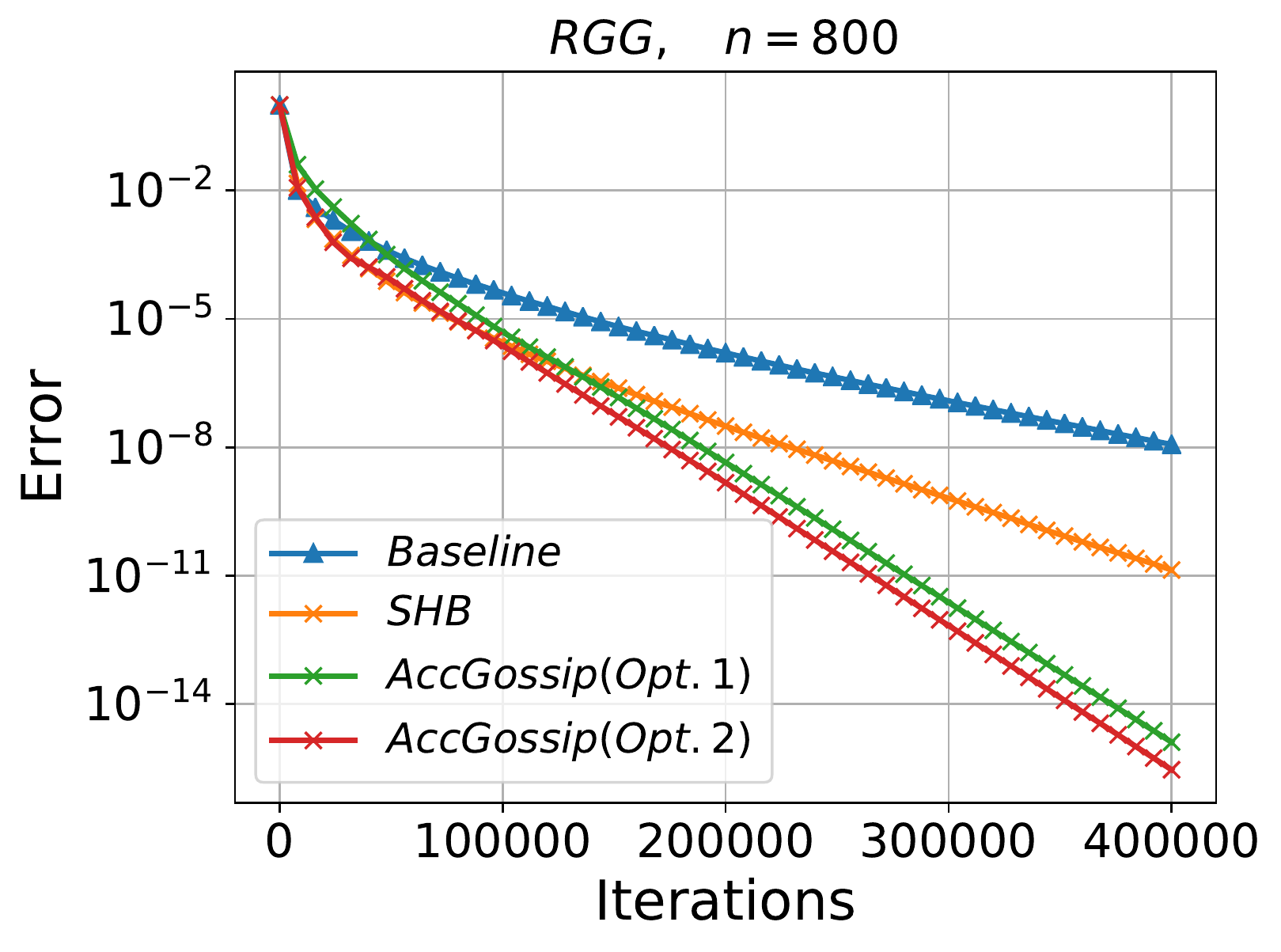}
\end{subfigure}\\
\begin{subfigure}{.24\textwidth}
  \centering
  \includegraphics[width=1\linewidth]{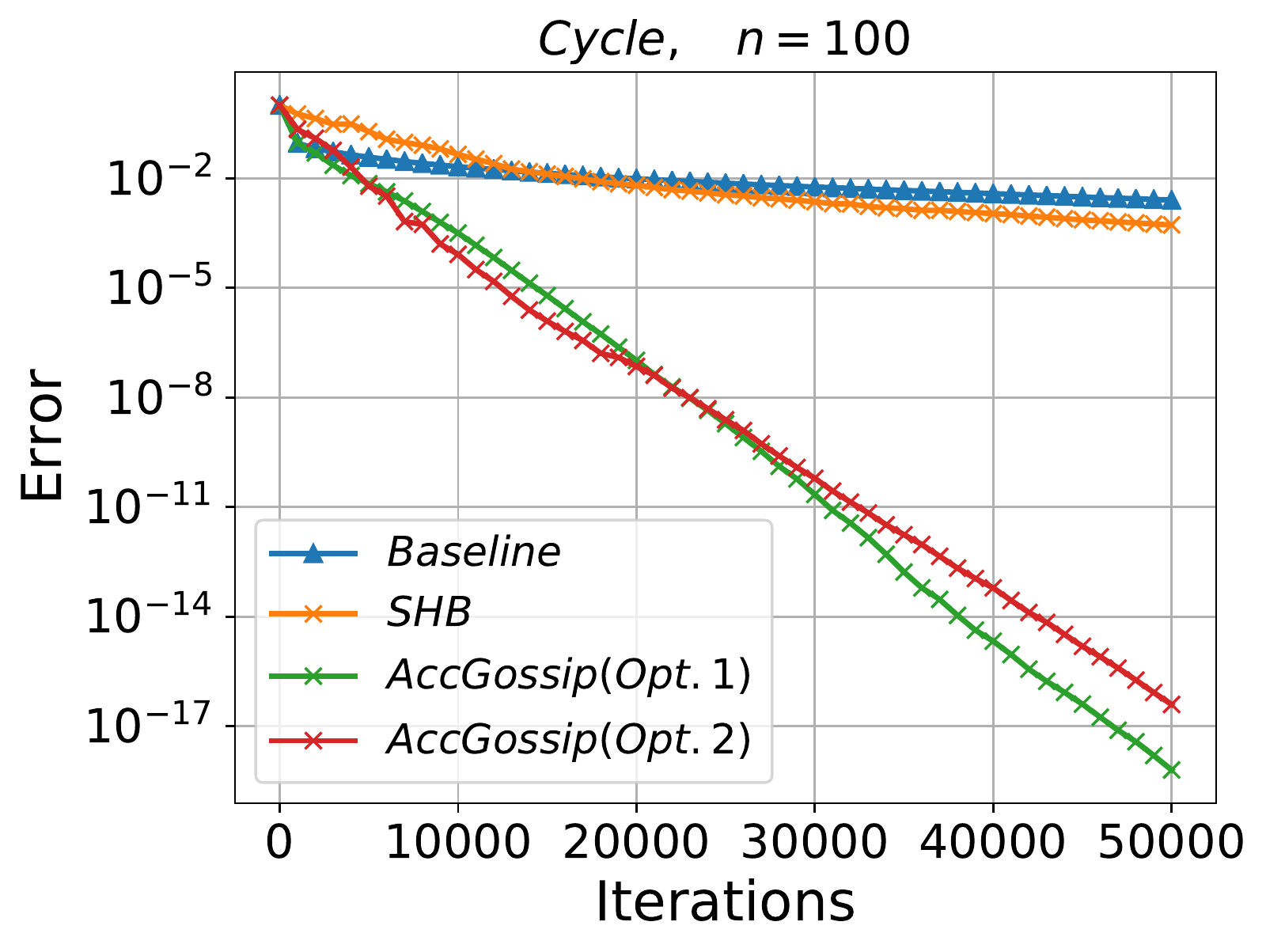}
\end{subfigure}%
\begin{subfigure}{.24\textwidth}
  \centering
  \includegraphics[width=1\linewidth]{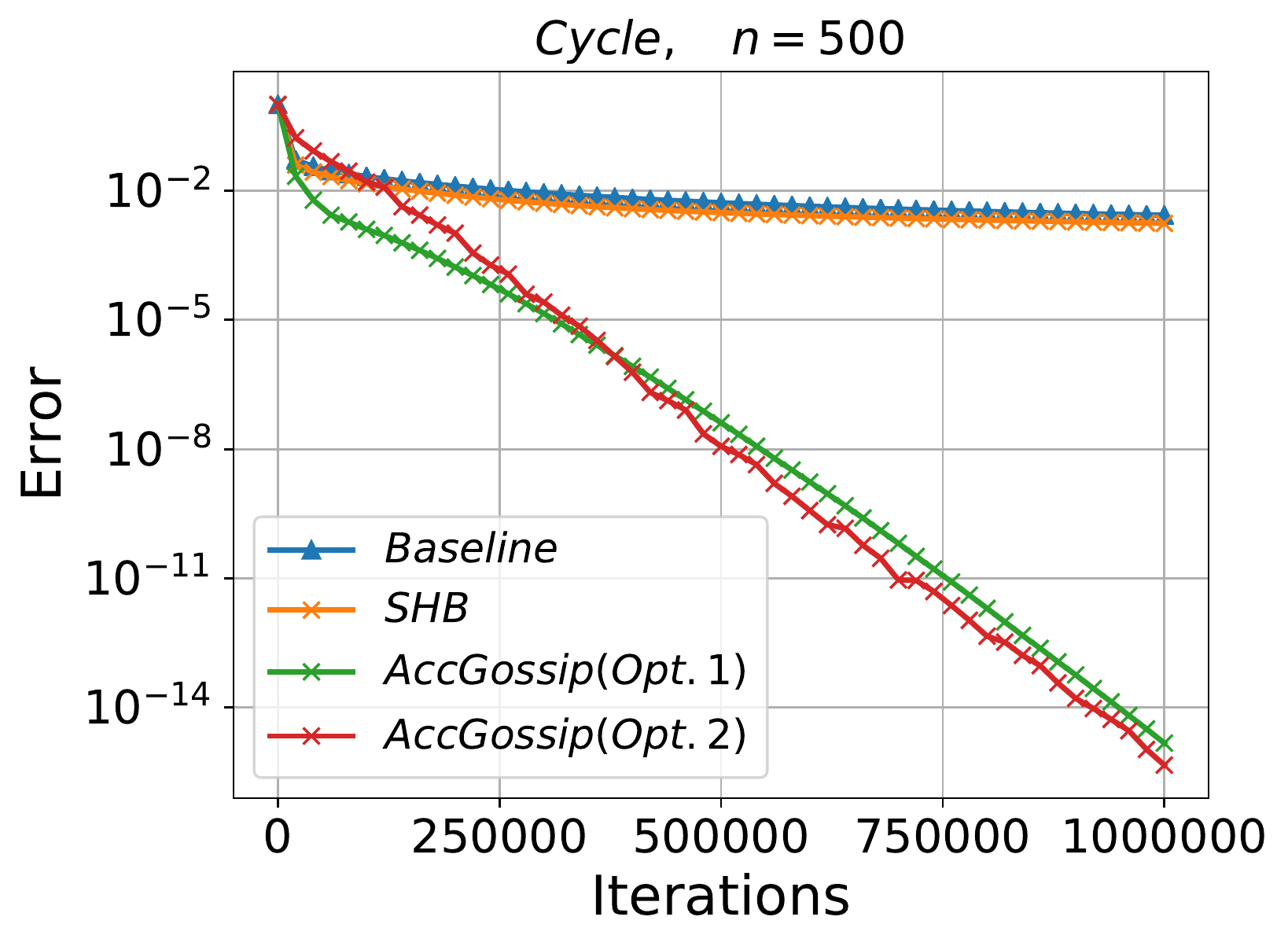}
\end{subfigure}
\caption{\footnotesize Performance of AccGossip in a 2-dimension grid, random geometric graph (RGG) and a cycle graph. The  Baseline method corresponds to the randomized pairwise gossip algorithm proposed in \cite{boyd2006randomized} and the SHB to the fast gossip algorithm proposed in \cite{loizou2018accelerated} ; The $n$ in the title of each plot indicates the number of nodes of the network. For the grid graph this is $n \times n$.}
\label{mRKomega1}
\end{figure}
\section{Conclusion and Future Research}
\label{sec:Conclusion}
We proposed novel provably accelerated randomized gossip algorithms for solving the AC problem. Our approach is based on connections established between the gossip algorithms and the Kaczmarz methods for solving linear systems. We believe that many novel and efficient gossip protocols can be discovered using results from the literature of Kaczmarz methods either by using different AC linear systems or using other Kaczmarz-type algorithms than the one presented in this manuscript. We speculate that the gossip algorithms presented in this work can be extended to the more general setting of minimizing the average of convex functions $(1/n) \sum_{i=1}^n f_i(x)$ in a decentralized way \cite{nedic2018network}. While preparing this work we become aware of \cite{hendrikx2018accelerated} where an accelerated gossip algorithm is developed for solving the dual of the best approximation problem \eqref{best approximation} using the accelerated coordinate descent method of \cite{nesterov2017efficiency}. A comparison of our protocols and the algorithm of \cite{hendrikx2018accelerated} is an ongoing research work. 

\vfill\pagebreak
\bibliographystyle{IEEEbib}
{\footnotesize\bibliography{AccGossip}}

\begin{thebibliography}{10}

\bibitem{degroot1974reaching}
Morris~H DeGroot,
\newblock ``Reaching a consensus,''
\newblock {\em Journal of the American Statistical Association}, vol. 69, no.
  345, pp. 118--121, 1974.

\bibitem{tsitsiklis1986distributed}
John Tsitsiklis, Dimitri Bertsekas, and Michael Athans,
\newblock ``Distributed asynchronous deterministic and stochastic gradient
  optimization algorithms,''
\newblock {\em IEEE transactions on automatic control}, vol. 31, no. 9, pp.
  803--812, 1986.

\bibitem{boyd2006randomized}
S.~Boyd, A.~Ghosh, B.~Prabhakar, and D.~Shah,
\newblock ``Randomized gossip algorithms,''
\newblock {\em IEEE Transactions on Information Theory}, vol. 14, no. SI, pp.
  2508--2530, 2006.

\bibitem{xiao2005scheme}
L.~Xiao, S.~Boyd, and S.~Lall,
\newblock ``A scheme for robust distributed sensor fusion based on average
  consensus,''
\newblock in {\em Information Processing in Sensor Networks, 2005. IPSN 2005.
  Fourth International Symposium on}. IEEE, 2005, pp. 63--70.

\bibitem{cybenko1989dynamic}
G.~Cybenko,
\newblock ``Dynamic load balancing for distributed memory multiprocessors,''
\newblock {\em J. Parallel Distrib. Comput.}, vol. 7, no. 2, pp. 279--301,
  1989.

\bibitem{freris2012fast}
N.M. Freris and A.~Zouzias,
\newblock ``Fast distributed smoothing of relative measurements,''
\newblock in {\em Decision and Control (CDC), 2012 IEEE 51st Annual Conference
  on}. IEEE, 2012, pp. 1411--1416.

\bibitem{dimakis2010gossip}
A.G. Dimakis, S.~Kar, J.M.F. Moura, M.G. Rabbat, and A.~Scaglione,
\newblock ``Gossip algorithms for distributed signal processing,''
\newblock {\em Proceedings of the IEEE}, vol. 98, no. 11, pp. 1847--1864, 2010.

\bibitem{zouzias2015randomized}
A.~Zouzias and N.M. Freris,
\newblock ``Randomized gossip algorithms for solving {Laplacian} systems,''
\newblock in {\em Control Conference (ECC), 2015 European}. IEEE, 2015, pp.
  1920--1925.

\bibitem{liu2013analysis}
J.~Liu, B.D.O. Anderson, M.~Cao, and A.S. Morse,
\newblock ``Analysis of accelerated gossip algorithms,''
\newblock {\em Automatica}, vol. 49, no. 4, pp. 873--883, 2013.

\bibitem{olshevsky2014linear}
A.~Olshevsky,
\newblock ``Linear time average consensus on fixed graphs and implications for
  decentralized optimization and multi-agent control,''
\newblock {\em arXiv preprint arXiv:1411.4186}, 2014.

\bibitem{LoizouRichtarik}
N.~Loizou and P.~Richt\'{a}rik,
\newblock ``A new perspective on randomized gossip algorithms,''
\newblock in {\em 4th IEEE Global Conference on Signal and Information
  Processing (GlobalSIP)}, 2016.

\bibitem{nedic2018network}
A.~Nedi{\'c}, A.~Olshevsky, and M.~G. Rabbat,
\newblock ``Network topology and communication-computation tradeoffs in
  decentralized optimization,''
\newblock {\em Proceedings of the IEEE}, vol. 106, no. 5, pp. 953--976, 2018.

\bibitem{aybat2017decentralized}
N.~S. Aybat and M.~G{\"u}rb{\"u}zbalaban,
\newblock ``Decentralized computation of effective resistances and acceleration
  of consensus algorithms,''
\newblock in {\em Signal and Information Processing (GlobalSIP), 2017 IEEE
  Global Conference on}. IEEE, 2017, pp. 538--542.

\bibitem{dimakis2008geographic}
A.G. Dimakis, A.D. Sarwate, and M.J. Wainwright,
\newblock ``Geographic gossip: {Efficient} averaging for sensor networks,''
\newblock {\em IEEE Trans. Signal Process.}, vol. 56, no. 3, pp. 1205--1216,
  2008.

\bibitem{aysal2009broadcast}
T.C. Aysal, M.E. Yildiz, A.D. Sarwate, and A.~Scaglione,
\newblock ``Broadcast gossip algorithms for consensus,''
\newblock {\em IEEE Trans. Signal Process.}, vol. 57, no. 7, pp. 2748--2761,
  2009.

\bibitem{olshevsky2009convergence}
A.~Olshevsky and J.N. Tsitsiklis,
\newblock ``Convergence speed in distributed consensus and averaging,''
\newblock {\em SIAM J. Control Optim.}, vol. 48, no. 1, pp. 33--55, 2009.

\bibitem{hanzely2017privacy}
F.~Hanzely, J.~Kone{\v{c}}n{\'y}, N.~Loizou, P.~Richt{\'a}rik, and
  D.~Grishchenko,
\newblock ``Privacy preserving randomized gossip algorithms,''
\newblock {\em arXiv preprint arXiv:1706.07636}, 2017.

\bibitem{loizou2018accelerated}
Nicolas Loizou and Peter Richt{\'a}rik,
\newblock ``Accelerated gossip via stochastic heavy ball method,''
\newblock {\em Allerton Conference on Communication, Control, and Computing,
  [arXiv preprint arXiv:1809.08657]}, 2018.

\bibitem{RK}
T.~Strohmer and R.~Vershynin,
\newblock ``A randomized {K}aczmarz algorithm with exponential convergence,''
\newblock {\em J. Fourier Anal. Appl.}, vol. 15, no. 2, pp. 262--278, 2009.

\bibitem{needell2010randomized}
D.~Needell,
\newblock ``Randomized {K}aczmarz solver for noisy linear systems,''
\newblock {\em BIT Numerical Mathematics}, vol. 50, no. 2, pp. 395--403, 2010.

\bibitem{RBK}
D.~Needell and J.A. Tropp,
\newblock ``Paved with good intentions: analysis of a randomized block
  {K}aczmarz method,''
\newblock {\em Linear Algebra and its Applications}, vol. 441, pp. 199--221,
  2014.

\bibitem{eldar2011acceleration}
Y.C. Eldar and D.~Needell,
\newblock ``Acceleration of randomized {K}aczmarz method via the
  {Johnson--Lindenstrauss} lemma,''
\newblock {\em Numerical Algorithms}, vol. 58, no. 2, pp. 163--177, 2011.

\bibitem{MaConvergence15}
A.~Ma, D.~Needell, and A.~Ramdas,
\newblock ``Convergence properties of the randomized extended {G}auss-{S}eidel
  and {K}aczmarz methods,''
\newblock {\em SIAM Journal on Matrix Analysis and Applications}, vol. 36, no.
  4, pp. 1590--1604, 2015.

\bibitem{zouzias2013randomized}
A.~Zouzias and N.M. Freris,
\newblock ``Randomized extended {K}aczmarz for solving least squares,''
\newblock {\em SIAM. J. Matrix Anal. \& Appl.}, vol. 34, no. 2, pp. 773--793,
  2013.

\bibitem{l2015randomized}
D.~Needell, R.~Zhao, and A.~Zouzias,
\newblock ``Randomized block {Kaczmarz} method with projection for solving
  least squares,''
\newblock {\em Linear Algebra and its Applications}, vol. 484, pp. 322--343,
  2015.

\bibitem{schopfer2016linear}
F.~Sch{\"o}pfer and D.A. Lorenz,
\newblock ``Linear convergence of the randomized sparse {K}aczmarz method,''
\newblock {\em arXiv preprint arXiv:1610.02889}, 2016.

\bibitem{liu2016accelerated}
J.~Liu and S.~Wright,
\newblock ``An accelerated randomized {Kaczmarz} algorithm,''
\newblock {\em Mathematics of Computation}, vol. 85, no. 297, pp. 153--178,
  2016.

\bibitem{loizou2017linearly}
N.~Loizou and P.~Richt{\'a}rik,
\newblock ``Linearly convergent stochastic heavy ball method for minimizing
  generalization error,''
\newblock {\em NIPS-Workshop on Optimization for Machine Learning [arXiv
  preprint arXiv:1710.10737]}, 2017.

\bibitem{gower2015randomized}
R.M. Gower and P.~Richt{\'a}rik,
\newblock ``Randomized iterative methods for linear systems,''
\newblock {\em SIAM. J. Matrix Anal. \& Appl.}, vol. 36, no. 4, pp. 1660--1690,
  2015.

\bibitem{gower2015stochastic}
R.M. Gower and P.~Richt{\'a}rik,
\newblock ``Stochastic dual ascent for solving linear systems,''
\newblock {\em arXiv preprint arXiv:1512.06890}, 2015.

\bibitem{loizou2017momentum}
N.~Loizou and P.~Richt{\'a}rik,
\newblock ``Momentum and stochastic momentum for stochastic gradient, newton,
  proximal point and subspace descent methods,''
\newblock {\em arXiv preprint arXiv:1712.09677}, 2017.

\bibitem{nesterov2012efficiency}
Y.~Nesterov,
\newblock ``Efficiency of coordinate descent methods on huge-scale optimization
  problems,''
\newblock {\em SIAM Journal on Optimization}, vol. 22, no. 2, pp. 341--362,
  2012.

\bibitem{tu2017breaking}
S.~Tu, S.~Venkataraman, A.C. Wilson, A.~Gittens, M.I. Jordan, and B.~Recht,
\newblock ``Breaking locality accelerates block {Gauss}-{Seidel},''
\newblock in {\em ICML}, 2017.

\bibitem{gower2018accelerated}
R.M Gower, Filip H., P.~Richt{\'a}rik, and S.~Stich,
\newblock ``Accelerated stochastic matrix inversion: general theory and
  speeding up bfgs rules for faster second-order optimization,''
\newblock {\em arXiv preprint arXiv:1802.04079}, 2018.

\bibitem{charalambous2016distributed}
T.~Charalambous, M.G. Rabbat, M.~Johansson, and C.N. Hadjicostis,
\newblock ``Distributed finite-time computation of digraph parameters: Left
  eigenvector, out-degree and spectrum,''
\newblock {\em IEEE Trans. Control of Network Systems}, vol. 3, no. 2, pp.
  137--148, June 2016.

\bibitem{cao2006accelerated}
M.~Cao, D.A. Spielman, and E.M. Yeh,
\newblock ``Accelerated gossip algorithms for distributed computation,''
\newblock in {\em Proc. of the 44th Annual Allerton Conference on
  Communication, Control, and Computation}, 2006, pp. 952--959.

\bibitem{hendrikx2018accelerated}
H.~Hendrikx, L.~Massouli{\'e}, and F.~Bach,
\newblock ``Accelerated decentralized optimization with local updates for
  smooth and strongly convex objectives,''
\newblock {\em arXiv preprint arXiv:1810.02660}, 2018.

\bibitem{nesterov2017efficiency}
Y.~Nesterov and S.U. Stich,
\newblock ``Efficiency of the accelerated coordinate descent method on
  structured optimization problems,''
\newblock {\em SIAM Journal on Optimization}, vol. 27, no. 1, pp. 110--123,
  2017.

\end{thebibliography}

\end{document}